\documentclass[11pt]{article}

\usepackage{sigsam, amsmath}

\issue{ACM 2010}

\articlehead{Extended abstract}

\titlehead{\textsf{simpcomp} -- A \textsf{GAP} toolbox for simplicial complexes}

\authorhead{Felix Effenberger and Jonathan Spreer}

\begin{document}

\title{ISSAC 2010 Software Presentation \\ \textsf{simpcomp} -- A \textsf{GAP} toolbox for simplicial complexes}

\author{Felix Effenberger and Jonathan Spreer \\
Fachbereich Mathematik\\
Institut f\"ur Geometrie und Topologie \\
Universit\"at Stuttgart, 70550 Stuttgart, Germany \\
\url{effenberger@mathematik.uni-stuttgart.de}, \\ \url{spreer@mathematik.uni-stuttgart.de}}

\date{}

\maketitle

\begin{abstract}
\textsf{simpcomp} \cite{simpcomp} is an extension (a so called \emph{package}) to \textsf{GAP} \cite{GAP4}, the well known system for computational discrete algebra. The package enables the user to compute numerous properties of (abstract) simplicial complexes, provides functions to construct new complexes from existing ones and an extensive library of triangulations of manifolds. For an introduction to the field of piecewise linear (PL) topology see the books \cite{Rourke72IntrPLTop} and \cite{Kuehnel95TightPolySubm}.

\medskip
The authors acknowledge support by the DFG: \textsf{simpcomp} has partly been developed within the DFG projects Ku 1203/5-2 and Ku 1203/5-3. 
\end{abstract}

\section{What is new}

\textsf{simpcomp} is a package for working with simplicial complexes in the \textsf{GAP} system. In contrast to the package \textsf{homology} \cite{Dumas04Homology} which focuses on simplicial homology computation, \textsf{simpcomp} claims to provide the user with a broader spectrum of functionality regarding simplicial constructions.   

\textsf{simpcomp} allows the user to interactively construct complexes and to compute their properties in the \textsf{GAP} shell. Furthermore, it makes use of \textsf{GAP}'s expertise in groups and group operations. For example,  automorphism groups and fundamental groups of complexes can be computed and  examined further within the \textsf{GAP} system. Apart from supplying a facet list, the user can as well construct simplicial complexes from a set of generators and a prescribed automorphism group -- the latter form being the common in which a complex is presented in a publication. This feature is to our knowledge unique to \textsf{simpcomp}. Furthermore, \textsf{simpcomp} as of Version 1.3.0 supports the construction of simplicial complexes of prescribed dimension, vertex number and transitive automorphism group as described in \cite{Lutz03TrigMnfFewVertVertTrans}, \cite{Casella01TrigK3MinNumVert}.

Furthermore, \textsf{simpcomp} has an extensive library of known triangulations of manifolds. This is the first time that they are easily accessible without having to look them up in the literature \cite{Kuehnel99CensusTight}, \cite{Casella01TrigK3MinNumVert}, or online \cite{Lutz08ManifoldPage}. This allows the user to work with many different known triangulations without having to construct them first. As of the current version 1.3.0 the library contains triangulations of roughly 650 manifolds and roughly 7000 pseudomanifolds, including all vertex transitive triangulations from \cite{Lutz08ManifoldPage}. Most properties that \textsf{simpcomp} can handle are precomputed for complexes in the library. Searching in the library is possible by the complexes' names as well as some of their properties (such as $f$-, $g$- and $h$-vectors and their homology).

\section{\textsf{simpcomp} benefits}

\textsf{simpcomp} is written entirely in the \textsf{GAP} scripting language, thus giving the user the possibility to see behind the scenes and to customize or alter \textsf{simpcomp} functions if needed. 

The main benefit when working with \textsf{simpcomp} over implementing the needed functions from scratch is that \textsf{simpcomp} encapsulates all methods and properties of a simplicial complex in a new \textsf{GAP} object type (as an abstract data type). This way, among other things, \textsf{simpcomp} can transparently cache properties already calculated, thus preventing unnecessary double calculations. It also takes care of the error-prone vertex labeling of a complex. 

\textsf{simpcomp} provides the user with functions to save and load the simplicial complexes to and from files and to import and export a complex in various formats (e.g. from and to \textsf{polymake/TOPAZ} \cite{Joswig00Polymake}, \textsf{Macaulay2} \cite{M2}, \LaTeX, etc.). 

In contrast to the software package \textsf{polymake} \cite{Joswig00Polymake} providing the most efficient algorithms for each task in form of a heterogeneous package (where algorithms are implemented in various languages), the primary goal when developing \textsf{simpcomp} was not efficiency (this is already limited by the \textsf{GAP} scripting language), but rather ease of use and ease of extensibility by the user in the \textsf{GAP} language with all its mathematical and algebraic capabilities. 

The package includes an extensive manual (see \cite{simpcomp}) in which all functionality of \textsf{simpcomp} is documented and also makes use of \textsf{GAP}'s built in help system such that all the documentation is available directly from the \textsf{GAP} prompt in an interactive way.  

\section{Some features that \textsf{simpcomp} supports}

\textsf{simpcomp} implements many standard and often needed functions for working with simplicial complexes. These functions can be roughly divided into three groups: (i) functions generating simplicial complexes (ii) functions to construct new complexes from old and (iii) functions calculating properties of  complexes -- for a full list of supported features see the documentation \cite{simpcomp}.

\textsf{simpcomp} furthermore implements a variety of functions connected to \emph{bistellar moves} (also known as \emph{Pachner moves} \cite{Pachner87KonstrMethKombHomeo}) on simplicial complexes. For example, \textsf{simpcomp} can be used to construct randomized spheres or randomize a given complex. Another prominent application of bistellar moves implemented in \textsf{simpcomp} is a heuristical algorithm that determines whether a simplicial complex is a \emph{combinatorial manifold} (i.e. that each link is PL homeomorphic to the boundary of the simplex). This algorithm was first presented by F.H.~Lutz and A.~Bj\"orner \cite{Bjoerner00SimplMnfBistellarFlips}. It uses a simulated annealing type strategy in order to minimize vertex numbers of triangulations while leaving the PL homeomorphism type invariant.

The package also supports \emph{slicings} of 3-manifolds (known as \emph{normal surfaces}, see \cite{Kneser29ClosedSurfIn3Mflds}, \cite{Haken61TheoNormFl}, \cite{Spreer10NormSurfsCombSlic}) and related constructions.

The current version of \textsf{simpcomp} is 1.3.0 (May 26th, 2010). On the roadmap for the next version 1.4.x which should appear still in 2010 are the support for \emph{simplicial blowups}, i.e. the resolutions of ordinary double points in combinatorial 4-pseudomanifolds. This functionality is to the authors' knowledge not provided by any other software package so far. Also, a closer interaction with the software system \textsf{Macaulay2} is planned.

\section{An example}

This section contains a small demonstration of the capabilities of \textsf{simpcomp} in form of an example construction.

M.~Casella and W.~Kühnel constructed a triangulated K3 surface with minimum number of 16 vertices in \cite{Casella01TrigK3MinNumVert}. They presented it in terms of the complex obtained by the automorphism group $G\cong AGL(1,\mathbb{F}_{16})$ given by the five generators
\begin{equation*}
	\scriptsize
	G=\left\langle
	\begin{array}{c}
		(1\,2)(3\,4)(5\,6)(7\,8)(9\,10)(11\,12)(13\,14)(15\,16),\\
		(1\,3)(2\,4)(5\,7)(6\,8)(9\,11)(10\,12)(13\,15)(14\,16),\\
		(1\,5)(2\,6)(3\,7)(4\,8)(9\,13)(10\,14)(11\,15)(12\,16),\\
		(1\,9)(2\,10)(3\,11)(4\,12)(5\,13)(6\,14)(7\,15)(8\,16),\\
		(2\,13\,15\,11\,14\,3\,5\,8\,16\,7\,4\,9\,10\,6\,12)
	\end{array}
	\right\rangle,
\end{equation*}
acting on the two generating simplices $\Delta_1=\langle 2,3,4,5,9 \rangle$ and $\Delta_2=\langle 2,5,7,10,11 \rangle$. It turned out to be a non-trivial problem 
to show (i) that the complex obtained is a combinatorial 4-manifold and (ii) to show that it is homeomorphic to a K3 surface as topological 4-manifold.

This turns out to be a rather easy task using \textsf{simpcomp}, as will be shown below. We will fire up \textsf{GAP}, load \textsf{simpcomp} and then construct the complex from its representation given above:

{\scriptsize
\begin{verbatim}
gap> LoadPackage("simpcomp");; #load the package
Loading simpcomp 1.1.21
by F.Effenberger and J.Spreer
http://www.igt.uni-stuttgart.de/LstDiffgeo/simpcomp
gap> SCInfoLevel(0);; #suppress simpcomp info messages 
gap> G:=Group((1,2)(3,4)(5,6)(7,8)(9,10)(11,12)(13,14)(15,16),
> (1,3)(2,4)(5,7)(6,8)(9,11)(10,12)(13,15)(14,16),
> (1,5)(2,6)(3,7)(4,8)(9,13)(10,14)(11,15)(12,16),
> (1,9)(2,10)(3,11)(4,12)(5,13)(6,14)(7,15)(8,16),
> (2,13,15,11,14,3,5,8,16,7,4,9,10,6,12));;
gap> K3:=SCFromGenerators(G,[[2,3,4,5,9],[2,5,7,10,11]]);
[SimplicialComplex

 Properties known: AutomorphismGroup, AutomorphismGroupSize, 
                   AutomorphismGroupStructure, AutomorphismGroupTransitivity, 
                   Dim, Facets, Generators, Name, VertexLabels.

 Name="complex from generators under group ((C2 x C2 x C2 x C2) : C5) : C3"
 Dim=4
 AutomorphismGroupSize=240
 AutomorphismGroupStructure="((C2 x C2 x C2 x C2) : C5) : C3"
 AutomorphismGroupTransitivity=2

/SimplicialComplex]
\end{verbatim}
}
\noindent
We first compute the $f$-vector, the Euler characteristic and the homology groups of \texttt{K3}:
{\scriptsize
\begin{verbatim}
gap> K3.F;
[ 16, 120, 560, 720, 288 ]
gap> K3.Chi;
24
gap> K3.Homology;
[ [ 0, [  ] ], [ 0, [  ] ], [ 22, [  ] ], [ 0, [  ] ], [ 1, [  ] ] ]
\end{verbatim}
}
\noindent
Now we verify that the complex \texttt{K3} is a combinatorial manifold using the heuristic algorithm based on bistellar moves described above:
{\scriptsize
\begin{verbatim}
gap> K3.IsManifold;
true
\end{verbatim}
}
\noindent
In a next step we compute the parity and the signature of the intersection form of the complex \texttt{K3}:
{\scriptsize
\begin{verbatim}
gap> K3.IntersectionFormParity;
0
gap> K3.IntersectionFormSignature;
[ 22, 3, 19 ]
\end{verbatim}
}
\noindent
This means that the intersection form of the complex \texttt{K3} is even. It has dimension $22$ and  signature $19-3=16$. Furthermore, \texttt{K3} is simply connected as can either be verified by showing that the fundamental group is trivial or by checking that the complex is $3$-neighborly:  
{\scriptsize
\begin{verbatim}
gap> K3.FundamentalGroup;
<fp group with 105 generators>
gap> Size(last);
1
gap> K3.Neighborliness;
3
\end{verbatim}
}
\noindent
It now follows from a theorem of M.~Freedman \cite{Freedman82Top4DimMnf} that the complex is in fact homeomorphic to a K3 surface because it has the same (even) intersection form. Furthermore, \texttt{K3} is a tight triangulation \cite{Kuehnel99CensusTight, Kuehnel95TightPolySubm} as can be verified as follows:
{\scriptsize
\begin{verbatim}
gap> K3.IsTight;
#I  SCIsTight: complex is (k+1)-neighborly 2k-manifold and thus tight.
true
\end{verbatim}
We can also have a look at the multiplicity vectors of the perfect polyhedral Morse function \cite[Sec.~3B]{Kuehnel95TightPolySubm} that orders the vertices $v_1,\dots,v_{16}$ linearly:
{\scriptsize
\begin{verbatim}
gap> SCMorseIsPerfect(K3,[1..16]);
true	
gap> SCMorseMultiplicityVector(K3,[1..16]);
[ [ 1, 0, 0, 0, 0 ], [ 0, 0, 0, 0, 0 ], [ 0, 0, 0, 0, 0 ], [ 0, 0, 1, 0, 0 ], 
  [ 0, 0, 2, 0, 0 ], [ 0, 0, 1, 0, 0 ], [ 0, 0, 4, 0, 0 ], [ 0, 0, 3, 0, 0 ], 
  [ 0, 0, 3, 0, 0 ], [ 0, 0, 4, 0, 0 ], [ 0, 0, 1, 0, 0 ], [ 0, 0, 2, 0, 0 ], 
  [ 0, 0, 1, 0, 0 ], [ 0, 0, 0, 0, 0 ], [ 0, 0, 0, 0, 0 ], [ 0, 0, 0, 0, 1 ] ]
\end{verbatim}
Finally, instead of constructing the triangulation \texttt{K3} from scratch, we could have also loaded it directly from the library, saving us a lot of work. We now load the complex from the library and verify that the complex from the library is combinatorially isomorphic to the complex \texttt{K3} that we constructed -- note below the two different searching methods provided by the library:
{\scriptsize
\begin{verbatim}
gap> SCLib.SearchByName("K3");      
[ [ 551, "K3 surface" ] ]
gap> SCLib.SearchByAttribute("Dim=4 and F[3]=Binomial(F[1],3)");
[ [ 11, "CP^2 (VT)" ], [ 551, "K3 surface" ] ]
gap> M:=SCLib.Load(551);;
gap> M.IsIsomorphic(K3);
true
\end{verbatim}

\normalsize
\end{document}